\documentclass[11pt]{article}
\usepackage{amsmath,amssymb,eucal}
\usepackage[dvips]{graphicx}

\newtheorem{thm}{\bf Theorem}[section]

\newtheorem{lemma}[thm]{\bf Lemma}

\newtheorem{rem}[thm]{\bf Remark}

\newtheorem{example}[thm]{\bf Example}
\newtheorem{defini}[thm]{\bf Definition}

\def\qed{$\Box$}

\begin{document}

\title{
Generalized rational blow-down, torus knots, \\
and Euclidean algorithm
}
\author{Yuichi YAMADA}
\date{\today}
%\date{ 2007}
\footnotetext[0]{%
}
\footnotetext[1]{%
This work was partially supported by Grant-in-Aid for
Scientific Research No.18740029, 
Japan Society for the Promotion of Science. 
} 
\maketitle
%%%%%%%%%%%%%%%%   
\begin{abstract}{
We construct a Kirby diagram of the rational homology ball
used in ^^ ^^ generalized rational blow-down" developed by Jongil Park.
The diagram consists of a dotted circle and a torus knot.
The link is simpler, but the parameters are a little complicate.
Euclidean Algorithm is used three times 
in the construction and the proof.
}\end{abstract}
%%%%%%%%%%%%%%%

%%%%%%%%%%%%%%%%%%%%%%%%%%%%%%%%%%%%%%%%%%
%%%%%%%%%%%%%%%%      section 1      %%%%%%%%%%%%%%%%%
%%%%%%%%%%%%%%%%%%%%%%%%%%%%%%%%%%%%%%%%%%
\section{Main theorem}\label{sec:intro}
%
%
%%%%%%%%%%%%%%%%%%%%%%%%%%%%%%%%%
\begin{figure}
\begin{center}
\includegraphics[scale=0.7]{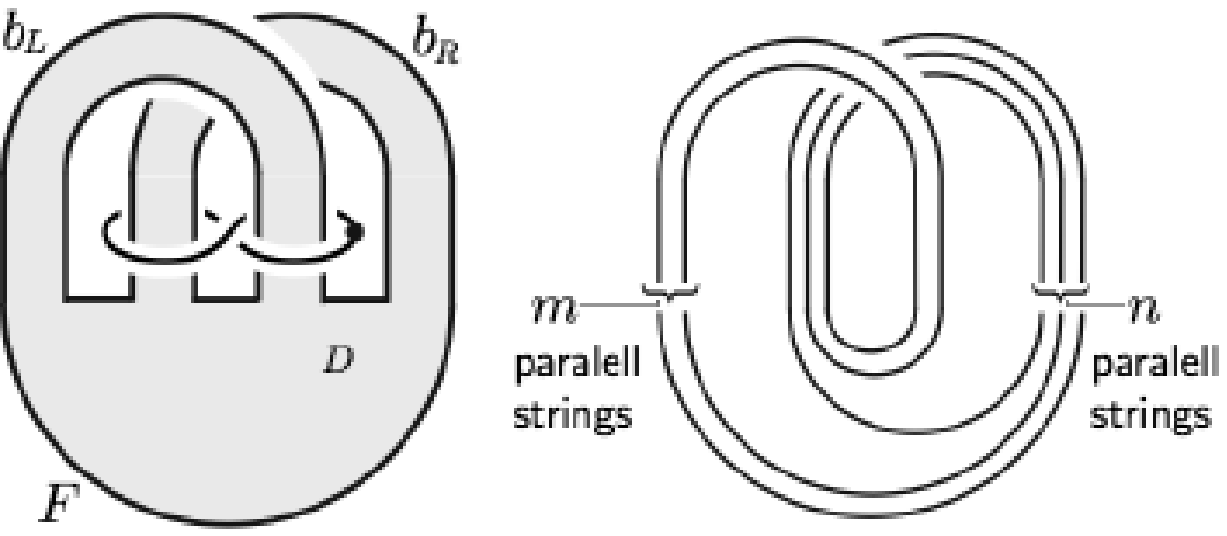} \quad
\includegraphics[scale=0.7]{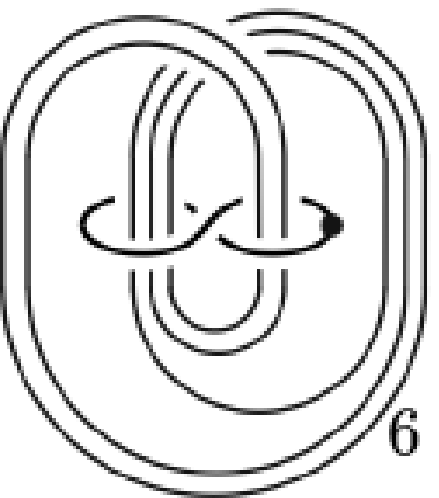} \\
\caption{\, $F$, \ $k(m,n)$ in $F$ 
\qquad ex. $k(2,3) \cup u$ ($= L_{5,2}$)}
\label{fig:knot}
\end{center}
\end{figure}
%%%%%%%%%%%%%%%%%%%%%%%%%%%%%%%%%
%
%
For a coprime pair $(m,n)$ of positive integers,
we take a simple closed curve $k(m,n)$
in the standardly embedded once-punctured torus $F$ in $S^3$
as in Figure~\ref{fig:knot}.
We study the Kirby diagram $k(m,n) \cup u$:
the component $k(m,n)$ is 
a torus knot $T(m,n)$ with $(mn)$-framing,
and $u$ is a dotted unknoted circle
(It is a 1-handle, see \cite{A,AK} and \cite[p.168]{GS})
in the complement of $F$.
This diagram defines a rational homology ball that has
cyclic fundamental group of order $(m+n)$.
It has a symmetry: $k(n,m) \cup u = k(m,n) \cup u$.
\par
In the next section, for a given coprime pair $(p,q)$ of
positive integers with $1 \leq q < p$, 
we will construct an involutive symmetric function $A$ by Algorithm,
to decide 
(another) coprime pair $(m, n) = A(p-q, q)$ satisfying $m + n =p$.
It holds that $A(p-1, 1) = (p-1, 1)$,
$A(p-2, 2) = ((p-1)/2, (p+1)/2)$ for odd $p$, 
see Lemma~\ref{lemma:A} and \ref{lemma:Amnst}.
Now, we let $L_{p,q}$ denote Kirby diagram $k(A(p-q, p)) \cup u$.
Our main theorem is:
%
%
%%%%%%%%%%%%%%%%%%%%%%%%%%%%%%
%%%%%%%%%%%%%%%%%%%%%%%%%%%%%%
\begin{thm}~\label{thm:main}
{\sl
For any coprime pair $(p,q)$ of positive integers
with $1 \leq q < p$, 
the boundary of the rational homology ball
described by $L_{p,q}$ defined above,
is a lens space $L(p^2, pq-1)$. 
}
\end{thm}
%%%%%%%%%%%%%%%%%%%%%%%%%%%%%%
%%%%%%%%%%%%%%%%%%%%%%%%%%%%%%
%
%
Thus, we can regard $L_{p,q}$ as a description of the
rational homology ball $B_{p,q}$ in 
{\it general rational blow-down} defined by J.~Park in \cite{P}
applying \cite{FS2} via \cite{CH}.
It is the operation ^^ ^^ cut out $C_{p,q}$ and paste $B_{p,q}$" on 
a $4$-manifold,
where $C_{p,q}$ is the negative definite 
plumbed $4$-manifold corresponding to the 
weighted graph in Figure~2.
The weights $(-c_i)$'s ($-c_i \leq -2$ for each $i$) are defined 
by the continued fraction expansion : 
$p^2/ (pq-1) = 
[c_0, c_1, \cdots, c_N]
$, where
\[
[x_1, x_2, \cdots, x_n] :=
x_1
- \cfrac{1}{x_2 
- \cfrac{1}{\ddots 
- \cfrac{1}{x_n}
}} \, .
\]
Thus $\partial C_{p,q} = L(p^2, pq-1) = \partial B_{p,q} $.
Note that $C_{p,q}$ has a symmetry: $C_{p,p-q} = C_{p,q}$,
corresponding to the {\it reverse} of the continued fraction
$[c_N, \cdots, c_1, c_0]$, 
and also to the homeomorphism
$L(p^2, p(p-q)-1) \cong L(p^2, pq-1)$. 
%
%
%%%%%%%%%%%%%%%%%%%%%%%%%%%%%%%%%
\begin{figure}
\begin{center}
\includegraphics[scale=0.3]{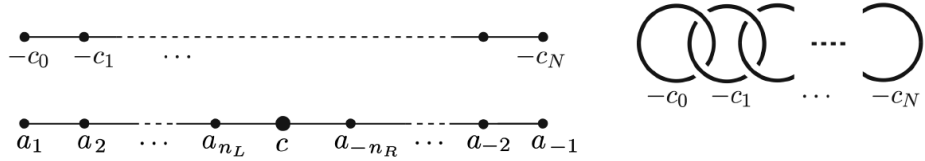} \quad
\caption{\, Plumbed manifold $C_{p,q}$}
\label{fig:plumb}
\end{center}
\end{figure}
%%%%%%%%%%%%%%%%%%%%%%%%%%%%%%%%%
%
%

Our strategy of the proof is:
First, in the next section, we will present Algorithm~A,
based on Euclidean algorithm of the pair $(p-q, q)$.
In the process, we construct 
a word $w(p-q,q)$ of $L$ and $R$, its ^^ ^^ reverse" $W(p-q,q)$,
decide integers $n_L, n_R$, and a finite sequence
\[
(a_1, a_2, \cdots, a_{nL},
\ c, \ 
a_{-nR}, \cdots, a_{-2}, a_{-1}
)
\] 
satisfying $a_i \leq -2$ (for each $i$) and $c \leq -4$.
Next in Section~3, we show that this sequence agrees to
$( -c_0, -c_1, \cdots, -c_N )$,
or its reverse $(-c_N, \cdots, -c_1, -c_0)$.
Finally, in Section~4, 
we prove Theorem~\ref{thm:main} by a sequence of Kirby calculus,
guided by the word $W(p-q,q)$ constructed in Algorithm A.
The process is related to the resolution (\cite{HKK, La}) 
of the singularity of the complex curve of type $z^m - w^n=0$,
or unknotting twisting sequence on torus knots, that is, 
Euclidean algorithm. \par

Note that $(mn)$-framed $T(m,n)$ is a kind of the most ^^ ^^ exceptional"
Dehn surgery, see \cite{M}.
Similar algorithm has been already discussed by the author
in \cite{Y3} (whose older version is in \cite{Y1}) 
in the study of exceptional Dehn surgery. 
The operation ^^ ^^ reverse" of the word at Step(2) in Algorithm~A,
is in contrast to the old results, and 
cause difficulty in the construction and the proof.
Some parts (ex. Figure~\ref{fig:seq})
are modification from the manuscript of \cite{Y3}, but
we rewrite them to make the present paper self-contained.
\par
%%%%
To the author's knowledge, 
descriptions of $B_{p,q}$ of some concrete $(p,q)$
can be seen, in \cite{SS} and \cite{R}. 
Our method is different from theirs: Non-trivial component
of the diagram of $B_{28,9}$ in \cite{SS} is 
$251$-framed $T(28,9)$ ($251 = 28 \cdot 9 -1$),
see also Remark~\ref{rem:final}.
%Note that Dehn surgery $T(a,b)$ with $(ab \pm 1)$-framing 
%is a lens space (\cite{M}).
Difference between theirs and ours looks like
a kind of ^^ ^^ dual", or ^^ ^^ complemental" in the sense 
that our function $A$ is involutive. \par
\medskip
The author would like to express sincere gratitude 
to Professor Jongil Park, for giving him the motivation of this research,
and some information on generalized rational blow-down, 
in and after one week lectures at The University of Tokyo, in June 2007.
The author would like to thank to Dr. Kouichi Yasui for 
valuable communication on rational blow-down and Kirby calculus.
\par
\medskip

%%%%%%%%%%%%%%%%%%%%%%%%%%%%%%%%%%%%%%%%%%
%%%%%%%%%%%%%%%%      section 2      %%%%%%%%%%%%%%%%%
%%%%%%%%%%%%%%%%%%%%%%%%%%%%%%%%%%%%%%%%%%
\section{Algorithm}\label{sec:algorithm}

Here we present the algorithm to 
define $(m,n) = A(p-q,q)$ via  
words $w(a,b)$, its reverse $W(a,b)$ (Here $a = p-q, b = q$),
define the integers $n_R$, $n_L$, and the sequence
\[
(a_1, a_2, \cdots, a_{nL},
\ c, \ 
a_{-nR}, \cdots, a_{-2}, a_{-1}
), \quad c = a_{n_L+1} + a_{-(n_R+1)} -2.
\]
This algorithm is closely related to the resolution (\cite{HKK, La}) 
of the singularity of the complex curve of type $z^a - w^b=0$,
that is, Euclidean algorithm.
We also show some formulas on $A$. 
\par
It may be curious, but we start with an example, which would help the readers.
%
%
%%%%%%%%%%%%%%%%%%%%%%%%%%%%%%
%%%%%%%%%%%%%%%%%%%%%%%%%%%%%%
\begin{example}~\label{ex:72}{\rm
$(a,b) = (7,2)$ (corresponding to $(p,q)=(9,2)$) \par 
\medskip
\begin{tabular}{lll}
$(a_i,b_i)$ : & $(7,2) \rightarrow_L (5,2) 
\rightarrow_L (3,2) 
\rightarrow_L (1,2) 
\rightarrow_R (1,1)$. 
& $w(7,2) = LLLR$. \\
$(m_i,n_i)$ : & $(1,1) \rightarrow_L (2,1) 
\rightarrow_L (3,1) 
\rightarrow_L (4,1) 
\rightarrow_R (4,5)$. 
& \ thus $A(7,2) = (4,5)$.\\
$(s_i,t_i)$ : & $(1,0) \rightarrow_L (1,0) 
\rightarrow_L (1,0) 
\rightarrow_L (1,0) 
\rightarrow_R (1,1)$. 
\end{tabular}
\quad 
\par \medskip
$n_L=3, n_R=1$, $W(7,2) = RLLL$.
\begin{center}
{\small
\begin{tabular}{l|ccccccccc}
$i$ & 
$a^{(i)}_{-2}$ & $a^{(i)}_{-1}$ & $a^{(i)}_0$ & $a^{(i)}_1$ & $a^{(i)}_2$
& $a^{(i)}_{3}$ & $a^{(i)}_{4}$  \\
\hline 
$0$ & & $-1$ & $-1$ & $-1$ & & & \\
$1$ & $-1$ & $-2$ & $-1$ & $-2$ & & &\\ 
$2$ & $-1$ & $-3$ & $-1$ & $-2$ & $-2$ & &\\ 
$3$ & $-1$ & $-4$ & $-1$ & $-2$ & $-2$ & $-2$ &\\ 
$4$ & $-1$ & $-5$ & $-1$ & $-2$ & $-2$ & $-2$ & $-2$ \\ 
\end{tabular}
}
{\small
\begin{tabular}{l|ccccccc}
$i$ & 
$\overline{a}^{(i)}_{-3}$ & $\overline{a}^{(i)}_{-2}$ & $\overline{a}^{(i)}_{-1}$ 
& $\overline{a}^{(i)}_0$ & $\overline{a}^{(i)}_1$ \\
\hline 
$0$ & & & & $-4$ & \\
$1$ & & & & $-5$ & $-2$ \\
$2$ & & & $-2$ & $-5$ & $-3$ \\
$3$ & & $-2$ & $-2$ & $-5$ & $-4$ \\
$4$ & $-2$ & $-2$ & $-2$ & $-5$ & $-5$ \\
\end{tabular}
}
\end{center}\par
We get the sequence $(-2, -2, -2, -5,-5)$ ($=C_{9,2}$),
and $[5,5,2,2,2] = 81/17$. See Figure~\ref{fig:bups}.
%For $(s_i,t_i)$, see the proof of Lemma~\ref{lemma:A}.
}
\end{example}
%%%%%%%%%%%%%%%%%%%%%%%%%%%%%%
%%%%%%%%%%%%%%%%%%%%%%%%%%%%%%
%
%
%%%%%%%%%%%%%%%%%%%%%%%%%%%%%%%%%
\begin{figure}[h]
\begin{center}
\includegraphics[scale=0.35]{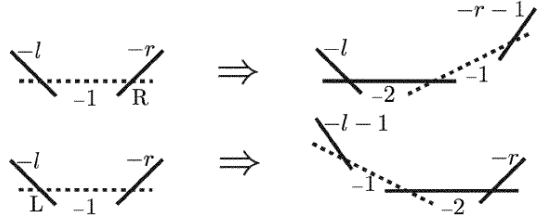} \quad
\caption{\, Blow-up at L or R}
\label{fig:LRbups}
\end{center}
\end{figure}
%%%%%%%%%%%%%%%%%%%%%%%%%%%%%%%%%
%
%
%%%%%%%%%%%%%%%%%%%%%%%%%%%%%%%%%
\begin{figure}[h]
\begin{center}
\includegraphics[scale=0.6]{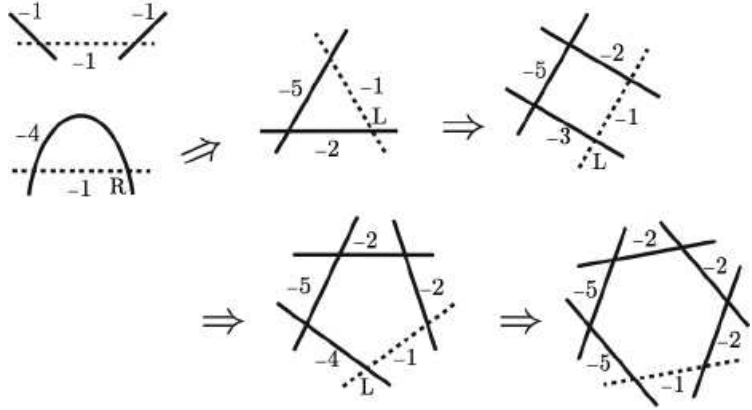} \quad
\caption{\, Blow-ups}
\label{fig:bups}
\end{center}
\end{figure}
%%%%%%%%%%%%%%%%%%%%%%%%%%%%%%%%%
%
%

\medskip \noindent
%%%%%%%%%%%%%%%%%%%%%%%%%%%%%%
{\bf Algorithm~A} 
\par \medskip \noindent
{\bf (1)} Euclidean algorithm: 
From the pair $(a,b)$ ($=(p-q,q)$), we
construct a word $w(a,b) = w_1 w_2 \cdots w_n$ of
two letters $L$(left) and $R$(right), and 
a sequence of the pair $\{ (m_i, n_i) \}$,
inductively, by the rule below: \par \noindent
Start with $(a_0,b_0) := (a,b)$, $(m_0,n_0) :=(1,1)$

\begin{center}
{\bf (LR Rule) }
\begin{tabular}{ll}
If $a_i > b_i$, 
& then $w_{i+1} := L$ and \\
& $(a_{i+1},b_{i+1}):=(a_i-b_i,b_i)$, \ 
$(m_{i+1}, n_{i+1}):= (m_i + n_i,n_i)$.\\[5pt]
If $a_i < b_i$,
& then $w_{i+1} := R$ and \\
& $(a_{i+1},b_{i+1}):=(a_i,b_i-a_i)$, \
$(m_{i+1}, n_{i+1}):= (m_i, n_i + m_i)$.
\end{tabular}
\end{center}
By coprime-ness of $(a,b)$, after some $N$ steps, the pair  
$(a_N,b_N)$ becomes to $(1,1)$, which is the end of this step.
%
%
%%%%%%%%%%%%%%%%%%%%%%%%%%%%%%
%%%%%%%%%%%%%%%%%%%%%%%%%%%%%%
\begin{defini}{\rm
We define $n_R$ (and $n_L$, respectively) 
as the number of $R$ (and $L$) in the word $w(a,b)$.
Thus $n_R +n_L = N$.
We define 
\[
A(a,b) := (m_N, n_N).
\]
}
\end{defini}
%%%%%%%%%%%%%%%%%%%%%%%%%%%%%%
%%%%%%%%%%%%%%%%%%%%%%%%%%%%%%
%
%
\par
\noindent
{\bf (2)} 
Let $W(a,b) = W_1W_2 \cdots W_N$ be the {\it reverse} 
of $w(a,b)$, i.e., $W_i = w_{N+1-i}$ for each $i$. It is easy to see
%
%
%%%%%%%%%%%%%%%%%%%%%%%%%%%%%%
%%%%%%%%%%%%%%%%%%%%%%%%%%%%%%
\begin{lemma}~\label{lemma:A}
{\sl 
Let $(a,b)$ a coprime pair of positive integers.
\begin{enumerate}
\item[(1)] If $A(a,b) =(m,n)$, then $W(a,b) = w(m,n)$, i.e., $W(a,b) =w(A(a,b))$.
\item[(2)] $A$ is involutive; If $A(a,b) = (m,n)$, then $A(m,n)=(a,b)$.
\item[(3)] $A$ is symmetric;  If $A(a,b) = (m,n)$, then $A(b,a) = (n,m)$.
\item[(4)] $A(a,1) = (a, 1)$. \ If $a$ is odd, 
$A(a,2) = (\frac{a+1}2, \frac{a+3}2)$.
\end{enumerate}
}
\end{lemma}
%%%%%%%%%%%%%%%%%%%%%%%%%%%%%%
%%%%%%%%%%%%%%%%%%%%%%%%%%%%%%
%
%
We go back to Algorithm~A. \par
\noindent
{\bf (3)} 
Next, starting with
\[
\{ a^{(0)}_{\ast} \} = 
(a^{(0)}_{-1}, a^{(0)}_{0}, a^{(0)}_{1}) := (-1, -1, -1),
\]
based on the blow-up diagram in Figure~\ref{fig:LRbups}
(see also Figure~\ref{fig:bups}),
we define the sequences $\{ a^{(i)}_{\ast} \}$ and $c^{(i)}$
($i=1,2, \cdots , N$) inductively:
For each $i$, $a^{(i)}_0 = -1$
and $c^{(i)} = a^{(i)}_{M(i)} + a^{(i)}_{m(i)} -2$, where
$M(i)$ (and $m(i)$, respectively)
is the maximum (or the minimum) in $\{ j \in {\bf Z} 
\vert \textrm{ $a^{(i)}_j$ is defined} \}$. Now, using $W_i$'s
(contrast to \cite{Y3}), 
\par \noindent
If $W_i=R$, then we define $\{ a^{(i)}_{\ast} \}$ as
\[
\begin{cases}
a^{(i)}_{j}:= a^{(i-1)}_{j} & 
\textrm{if \ $1 < j \leq M(i-1)$} \\
%\textrm{if $j>1$ and $a^{(i-1)}_{j}$ is defined,} \\
a^{(i)}_{1}:= a^{(i-1)}_{1}-1, \\
a^{(i)}_{-1}:= -2, \\
a^{(i)}_{j}:= a^{(i-1)}_{j+1} & 
\textrm{if \ $m(i-1)-1 \leq j <-1$}
%\textrm{if $j<-1$ and $a^{(i-1)}_{j+1}$ is defined.}
\end{cases},
\]
\noindent
If $W_i=L$, then we define $\{ a^{(i)}_{\ast} \}$ as
\[
\begin{cases}
a^{(i)}_{j}:= a^{(i-1)}_{j} & 
\textrm{if \ $m(i-1) \leq j < -1$} \\
%\textrm{if $j<-1$ and $a^{(i-1)}_{j}$ is defined,}\\
a^{(i)}_{-1}:= a^{(i-1)}_{-1}-1, \\
a^{(i)}_{1}:= -2, \\
a^{(i)}_{j}:= a^{(i-1)}_{j-1} & 
\textrm{if \ $1 < j \leq M(i-1)+1$} \\
%\textrm{if $j>1$ and $a^{(i-1)}_{j-1}$ is defined.}
\end{cases}.
\]
\noindent
{\bf (4)} 
For each integer $j$ with $-n_R \leq j \leq n_L$,
we define $a_j$ as $a^{(N)}_{j}$ 
in the sequence $\{ a^{(N)}_{\ast} \}$ 
obtained after the $N$-th step, where $N$ is the length 
of the word $W(a,b)$.
We also define $c := c^{(N)} = a_{n_L+1} + a_{-(n_R+1)} -2$. 
This is the end of Algorithm~A \hfill \qed 
\par \medskip
Related to Seifert fibration of $S^3$ whose regular fiber
is the torus knot $T(m,n)$, it is well-known:
\par
%
%
%%%%%%%%%%%%%%%%%%%%%%%%%%%%%%
%%%%%%%%%%%%%%%%%%%%%%%%%%%%%%
\begin{lemma}~\label{lemma:tknot}
{\sl
If $m < n$ (or  if $m > n$, respectively),
then $a_{-(n_R+1)} = -1$ (or $a_{n_L+1} = -1$), 
and 
\begin{eqnarray*}
& \left [
\vert a_{-(n_R+1)} \vert, 
\vert a_{-n_R} \vert,
\cdots,
\vert a_{-2} \vert,
\vert a_{-1} \vert
\right ]
& = 
\begin{cases}
m/n & \textrm{ if $m < n$} \\
n/m & \textrm{ if $m > n$}
\end{cases},
\\[10pt]
%%%
& \left [
\vert a_{(n_L+1)} \vert, 
\vert a_{n_L} \vert,
\cdots,
\vert a_{2} \vert,
\vert a_{1} \vert
\right ]
\hfill 
& = 
\begin{cases}
n/m & \textrm{ if $m < n$} \\
m/n & \textrm{ if $m > n$}
\end{cases}.
\end{eqnarray*}
}
\end{lemma}
%%%%%%%%%%%%%%%%%%%%%%%%%%%%%%
%%%%%%%%%%%%%%%%%%%%%%%%%%%%%%
%
%
Here, we add two formulas on $A$.
%
%
%%%%%%%%%%%%%%%%%%%%%%%%%%%%%%
%%%%%%%%%%%%%%%%%%%%%%%%%%%%%%
\begin{lemma}~\label{lemma:Amnst}
{\sl 
Let $(a,b)$ a coprime pair of positive integers.
Suppose $A(a,b) = (m,n)$, then 
\begin{enumerate}
\item[(1)] $m +n = a + b$ \par \noindent
\item[(2)] Let $s, t$ be the unique positive integers
that satisfies $mt - ns = -1$ and $0< s, t < a+b$,
then $s+t = b$.
\end{enumerate}
}
\end{lemma}
%%%%%%%%%%%%%%%%%%%%%%%%%%%%%%
%%%%%%%%%%%%%%%%%%%%%%%%%%%%%%
%
%
{\it Proof.} We go back to (LR Rule) in the construction of the function $A$.
(1) If $w_i = L$, then we took 
$(a_{i+1}, b_{i+1}) = (a_{i} - b_{i}, b_{i})$ and
$(m_{i+1}, n_{i+1}) = (m_{i} + n_{i}, n_{i})$. 
The equality 
$a_{i} n_{i} + b_{i} m_{i} = a+b$ is kept, in the process. 
It is kept also in the case $w_i =R$.
\par
(2) Similarly to $(m_i,n_i)$, 
we define $(s_i, t_i)$ inductively as: Starting $(s_0, t_0) = (1, 0)$,
if $w_i = L$ (or $R$, respectively), then we take
$(s_{i+1}, t_{i+1}) = (s_i + t_i, t_i)$
(or $(s_{i+1}, t_{i+1}) = (s_i, t_i + s_i )$).
Then the equality 
$m_{i} t_{i} - n_{i} s_{i} = -1$ is kept. 
We have $(s,t) = (s_N, t_N)$.
The equality 
$a_{i} t_{i} + b_{i} s_{i} = b$ is also kept.
\qed
\par \medskip
%%%%%%%%%%%%%%%%%%%%%%%%%
\noindent
{\bf (3)'} In addition, we present another algorithm to construct a
sequence $\{ \overline{a}_{\ast} \}$.
The author has been informed by J.~Park that 
%this is in the theory of ^^ ^^ the {\bf Q} Gorenstein singularity of class T".
this is a resolution graph of a quotient singularity of class T.
The resulting sequence $\{ \overline{a}_{\ast} \}$
will be agree to the sequence we have constructed in Algorithm~A.
\par
Starting with
\[
\{ \overline{a}^{(0)}_{\ast} \} = 
(\overline{a}^{(0)}_{0}) := (-4),
\]
we define the sequences $\{ \overline{a}^{(i)}_{\ast} \}$
($i=1,2, \cdots , n$) inductively.
We set
$M(i) := \max \{ j \in {\bf Z} 
\vert \textrm{ $\overline{a}^{(i)}_j$ is defined} \}$,
and $m(i) := \min \{ j \in {\bf Z} 
\vert \textrm{ $\overline{a}^{(i)}_j$ is defined} \}$.
\par \noindent
If $W_i=R$, then we define $\{ a^{(i)}_{\ast} \}$ as
\[
\begin{cases}
\overline{a}^{(i)}_{j+1}:= -2 & 
\textrm{ if  $j = M(i-1)$ } \\
\overline{a}^{(i)}_{j} \ \, := \overline{a}^{(i-1)}_{j} -1 & 
\textrm{ if  $j = m(i-1)$ } \\
\overline{a}^{(i)}_{j} \ \, := \overline{a}^{(i-1)}_{j} & 
\textrm{ if  $m(i-1) < j \leq M(i-1)$ } \\
\end{cases} ,
\]
If $W_i=L$, then we define $\{ a^{(i)}_{\ast} \}$ as
\[
\begin{cases}
\overline{a}^{(i)}_{j} \ \, := \overline{a}^{(i-1)}_{j} -1 & 
\textrm{ if  $j = M(i-1)$ } \\
\overline{a}^{(i)}_{j-1}:= -2 & 
\textrm{ if  $j = m(i-1)$ } \\
\overline{a}^{(i)}_{j} \ \, := \overline{a}^{(i-1)}_{j} & 
\textrm{ if  $m(i-1) \leq j < M(i-1)$ } \\
\end{cases}.
\]
Finally, we define $\overline{a}_j$ as $\overline{a}^{(N)}_{j}$ 
in the sequence $\{ \overline{a}^{(N)}_{\ast} \}$ 
obtained after the $N$-th step.
Note that $\overline{a}^{(i)}_0 = c^{(i)}$ and 
$\overline{a}_0 = c (= a_{n_L+1} + a_{-(n_R+1)} -2)$. 
\par
%
%
%%%%%%%%%%%%%%%%%%%%%%%%%%%%%%
%%%%%%%%%%%%%%%%%%%%%%%%%%%%%%
\begin{lemma}~\label{lemma:2seq}
{\sl
Two resulting sequence 
agrees to each other, i.e.,
\[
(a_1, a_2, \cdots, a_{n_L}, c , 
a_{-n_R}, \cdots, a_{-2}, a_{-1})
\ = \
(\overline{a}_{-n_L}, \cdots, \overline{a}_0, \cdots,  \overline{a}_{n_R}).
\]
}
\end{lemma}
%%%%%%%%%%%%%%%%%%%%%%%%%%%%%%
%%%%%%%%%%%%%%%%%%%%%%%%%%%%%%
%
%
{\it Proof.} (See Figure~\ref{fig:bups} again.) 
As a cyclic diagram,
$( a_1, a_2, \cdots, a_{nL}, c , $ 
$a_{-nR}, \cdots, a_{-2}, a_{-1})$ connected by $a_0 (= -1)$
in natural order, agrees to that of 
$(\overline{a}_{-n_L}, \cdots, \overline{a}_{n_R})$ 
connected by  
$\overline{a}_{-(n_L+1)} := -1 =: \overline{a}_{n_R+1}$.
We have the lemma. \qed
%
%
%%%%%%%%%%%%%%%%%%%%%%%%%%%%%%%%%
\par
\medskip

%%%%%%%%%%%%%%%%%%%%%%%%%%%%%%%%%%%%%%%%%%
%%%%%%%%%%%%%%%%      section 3      %%%%%%%%%%%%%%%%%
%%%%%%%%%%%%%%%%%%%%%%%%%%%%%%%%%%%%%%%%%%
\section{Sequence and lens space}\label{sec:sequence}

From a given coprime pair $(p,q)$, the sequence 
$(a_1, a_2, \cdots, a_{n_L}, c , 
a_{-n_R}, \cdots, a_{-2}, a_{-1})$
has been constructed in the previous section. Now, we show
%
%
%%%%%%%%%%%%%%%%%%%%%%%%%%%%%%
%%%%%%%%%%%%%%%%%%%%%%%%%%%%%%
\begin{lemma}~\label{lemma:Cfrac}
{\sl
The plumbed $4$-manifold of the weighted tree of 
the above sequence is diffeomorphic to $C_{p,q}$.
In other words,
the continued fraction expansion of 
$p^2/(pq-1)$ agrees to the sequence of the 
absolute values, up to reverse, i.e., it holds
\[
[ \vert a_1 \vert, \cdots, \vert a_{n_L} \vert , \vert c \vert, 
\vert a_{-n_R} \vert , \cdots, \vert a_{-1} \vert ]
\ = \
\dfrac{p^2}{pq-1},
\]
or
\[
[ \vert a_{-1} \vert, \cdots, \vert a_{-n_R} \vert , \vert c \vert, 
\vert a_{n_L} \vert , \cdots, \vert a_1 \vert ]
\ = \
\dfrac{p^2}{pq-1}.
\]
}
\end{lemma}
%%%%%%%%%%%%%%%%%%%%%%%%%%%%%%
%%%%%%%%%%%%%%%%%%%%%%%%%%%%%%
%
%

\begin{figure}
\begin{center}
\includegraphics[scale=0.35]{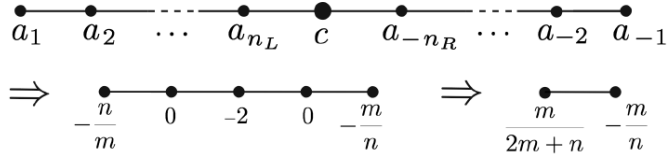} \quad
\caption{\, Lens space ($C_{p,q}$)}
\label{fig:lens}
\end{center}
\end{figure}

{\it Proof.} 
We also use $(m,n) =A(p-q,q)$, and $(s,t)$ satisfying $mt-ns=-1$
defined in Section~\ref{sec:algorithm}.
By Lemma~\ref{lemma:tknot}, up to reverse, we can contract 
the weighted graph as in Figure~\ref{fig:lens}.
%
%
%%%%%%%%%%%%%%%%%%%%%%%%%%%%%%%%%
%%%%%%%%%%%%%%%%%%%%%%%%%%%%%%%%%
%
%
In general, if two fractions at the vertices are 
$- \dfrac{\alpha_1}{\beta_1}$ and $- \dfrac{\alpha_2}{\beta_2}$,
the corresponding lens space is $L(P,Q)$ with
$P =  \alpha_1 \alpha_2 - \beta_1 \beta_2, \
Q \equiv  \alpha_1 \gamma_2 - \beta_1 \delta_2$ (mod $P$), 
where $\gamma_2, \delta_2$ are integers satisfying
$\alpha_2 \delta_2 - \beta_2 \gamma_2 = -1$.
In our case $(\alpha_2, \beta_2) = (m,n)$,
by Lemma~\ref{lemma:Amnst}(2), we set
$(\gamma_2, \delta_2) = (s, t)$. Thus,
\begin{eqnarray*}
P & = & m \cdot m + n \cdot (2m+n) = (m+n)^2, \\
Q & = & m \cdot s + (2m+n) \cdot t \\
& = & m(s+t) + mt + nt \ = m(s+t) + (ns -1) + nt\\
& = & (m+n)(s + t)-1.
\end{eqnarray*}
By Lemma~\ref{lemma:Amnst}(1) $m+n = a+ b =p$, and 
(2) $s+t = b =q$,
we have 
$P = p^2, Q= pq-1$. 
By the uniqueness of the continued fraction expansion
(with $a_i < -2$, $c<-2$, and $c_i <-2$),
we have the lemma. \qed
\par
\medskip

%%%%%%%%%%%%%%%%%%%%%%%%%%%%%%%%%%%%%%%%%%
%%%%%%%%%%%%%%%%      section 4      %%%%%%%%%%%%%%%%%
%%%%%%%%%%%%%%%%%%%%%%%%%%%%%%%%%%%%%%%%%%
\section{Proof of the Main Theorem}\label{sec:proof}

Let $(a,b)=(p-q,q)$ and $A(a,b) = A(p-q,q)=(m,n)$ as before.
Here we prove that the boundary of the 
rational ball described by Kirby diagram 
$L_{p,q} = k(m,n) \cup u$ is homeomorphic to $L(p^2, pq-1)$.\par

We have defined $F$ as a standardly embedded once-punctured 
torus in $S^3$, see Figure~1 again. 
It consists of a disk $D$ and two bands $b_L$ and $b_R$.
We took a simple closed curve 
$k(m,n)$ in $F$ as in Figure~1. 
The framing of $k(m,n)$ 
defined by the surface $F$ is $(mn)$.
From now on, we call such a framing {\it $F$-framing} 
(^^ ^^ surface framing"). 
%%%

\begin{figure}
\begin{center}
\includegraphics[scale=0.7]{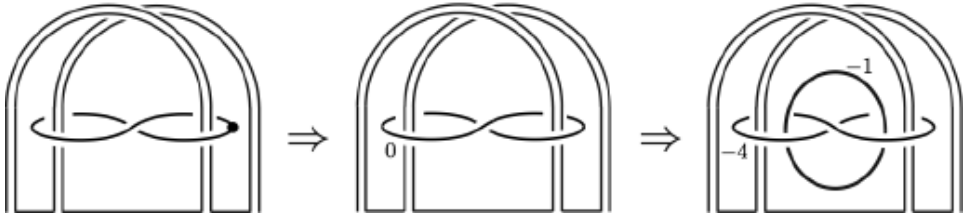} \quad
\caption{\, First Kirby move }
\label{fig:1Kmv}
\end{center}
\end{figure}

Our first Kirby move is in Figure~\ref{fig:1Kmv},
where, and from now on, we draw neither $D$ nor 
the components $k(m, n)$'s:
(1) Exchange the dotted circle $u$
to a $0$-framed same component, say $u'$.
This operation corresponds to a surgery
(cut out $S^1 \times D^3$ and paste $D^2 \times S^2$) 
in the interior of the rational ball, see \cite[p.7]{K} or \cite[p.168]{GS}.
Thus the boundary is unchanged.
(2) Blow-up. The central crossing is changed.
%
%
%%%%%%%%%%%%%%%%%%%%%%%%%%%%%%%%%
%%%%%%%%%%%%%%%%%%%%%%%%%%%%%%%%%
%
%

Before starting the next step, we define a notation:
$(\overline{m}_i, \overline{n}_i) := (m_{N+1-i}, n_{N+1-i})$,
where $\{ (m_j ,n_j) \}$ is the sequence of the pair
constructed in Step(1) in Algorithm~A. 
Thus, $(\overline{m}_0, \overline{n}_0) =(m,n)$ decreases to
$(\overline{m}_N, \overline{n}_N) =(1,1)$ guided by 
$w(m,n) = W(a,b) = W_1W_2 \cdots W_N$
(Lemma~{\ref{lemma:A}(1)), i.e., \underline{it holds} that \par
%%%%%%%%%%%%
\begin{center}
\begin{tabular}{ll}
If $\overline{m}_i < \overline{n}_i$, &
then $W_{i+1} = R$ and 
$(\overline{m}_{i+1},\overline{n}_{i+1}) 
= (\overline{m}_i, \overline{n}_i -\overline{m}_i)$, and \\
If $\overline{m}_i > \overline{n}_i$, &
then $W_{i+1} = L$ and 
$(\overline{m}_{i+1},\overline{n}_{i+1}) 
= (\overline{m}_i-\overline{n}_i, \overline{n}_i)$. 
\end{tabular}
\end{center}
%%%%%%%%%%%%

Next, guided by $W_1W_2 \cdots W_N$,
we move $F$ and the curve $k(m,n) = k(\overline{m}_0, \overline{n}_0)$
simultaneously in the total space $S^3$,
inductively ($i = 0,1,2, \cdots, N$):
If $W_{i+1} = R$ (i.e., $\overline{m}_i < \overline{n}_i$),
we move the \underline{left} band $b_L$ 
over the central $(-1)$-component
and slide over $b_R$ as in Figure~\ref{fig:seq}.
%
%
%%%%%%%%%%%%%%%%%%%%%%%%%%%%%%%%%
\begin{figure}
\begin{center}
\includegraphics[scale=1]{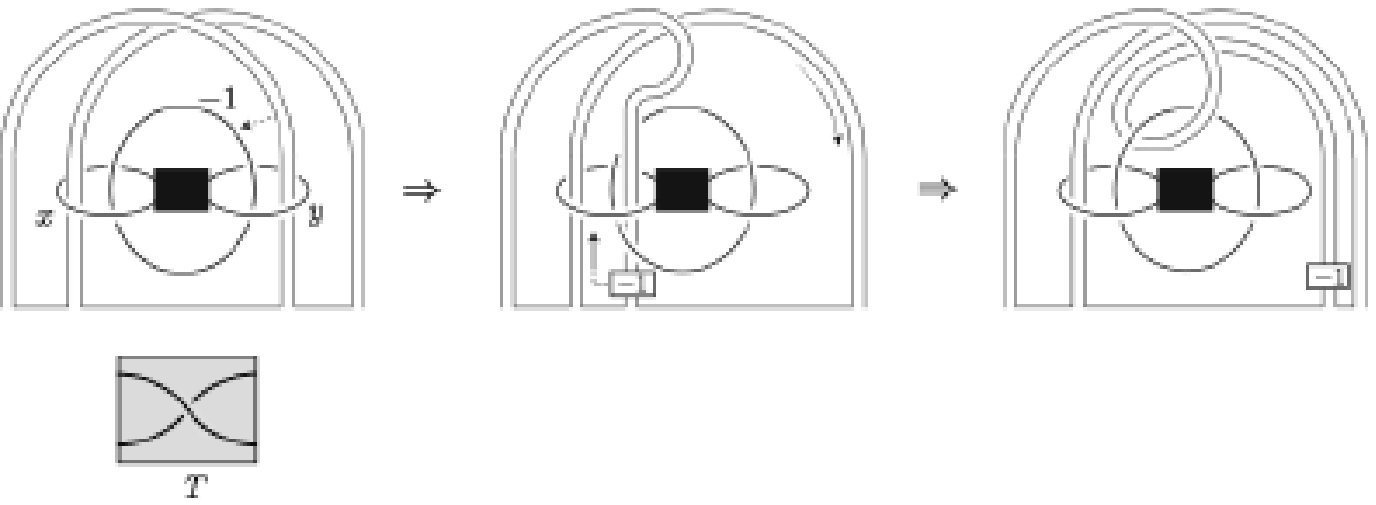} \\
\includegraphics[scale=1]{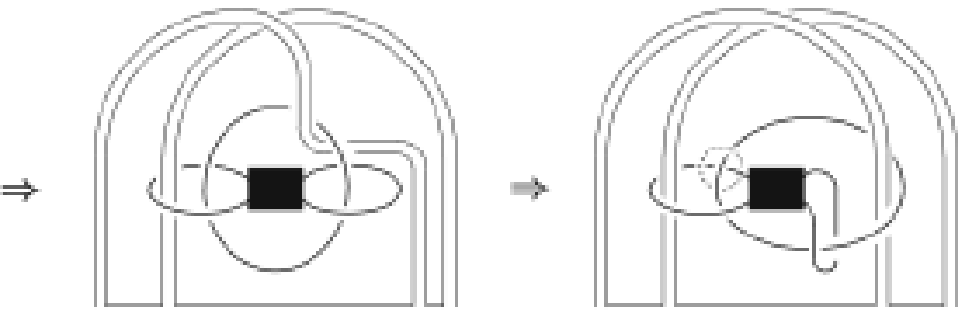} \\
\includegraphics[scale=1]{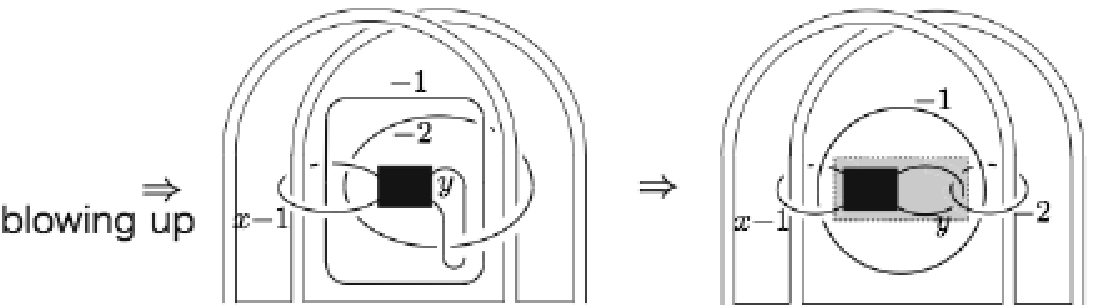} \\
\caption{\, Operation (R case)}
\label{fig:seq}
\end{center}
\end{figure}
%%%%%%%%%%%%%%%%%%%%%%%%%%%%%%%%%
%
%
In the black box, in the first step,
we take a tangle $T$ (two sub-arcs of the same component 
$(-4)$-framed $u'$), and
in the second or later steps,
we take the tangle 
that appeared in the gray box
at the end of the previous step, inductively.
In the case $W_{i+1} =L$, exchange the right and the left, 
but it is similar by symmetry.
Note that after a set of operation in Figure~\ref{fig:seq}, which includes
one blow-up, $F$ comes back to the starting position and 
$k(\overline{m}_i, \overline{n}_i)$ is changed 
to $k(\overline{m}_i, \overline{n}_i - \overline{m}_i)$ in $R$ case or 
to $k(\overline{m}_i - \overline{n}_i, \overline{n}_i)$ in $L$ case,
that is, to $k(\overline{m}_{i+1}, \overline{n}_{i+1})$ in either case
and new $(-1)$-component appears for the next step. 
Note that the relation 
^^ ^^ $F$-framing of $k(\overline{m}_i, \overline{n}_i)$
is $(\overline{m}_i \overline{n}_i)$"
is kept during the process.\par

After $N$ steps
($N$ is the length of the word $w(m,n)$, 
see also Lemma~\ref{lemma:A}(1)),
the diagram that we hoped appears at the black box, 
because this sequence of blow-ups exactly
same with the construction (3)' of $\{ \overline{a}_* \}$
in Section~\ref{sec:algorithm}.
By Lemma~\ref{lemma:2seq} and \ref{lemma:Cfrac},
it is the diagram of $C_{p,q}$, up to reverse.

The final $(-1)$-curve 
$\gamma$ and a $(+1)$-framed curve 
$\gamma' :=  k(\overline{m}_N, \overline{n}_N) = k(1,1)$ in $F$.
Sliding $\gamma'$ over $\gamma$, 
we can cancel them. 
The proof of Theorem~\ref{thm:main} is completed.  \qed
\par \medskip

Our proof shows:
%
%
%%%%%%%%%%%%%%%%%%%%%%%%%%%%%%
%%%%%%%%%%%%%%%%%%%%%%%%%%%%%%
\begin{rem}~\label{rem:final}
{\rm
By the Kirby calculus in our proof,
the $0$-framed meridian of $k(m,n)$ in $L_{p,q}$ comes
to $(-1)$-framed $\gamma$ in the diagram of $C_{p,q}$ 
in Figure~\ref{fig:total}.
This is (as a link component, at least) different from 
the example in \cite{SS},
or the one obtained by the method in \cite{R} and \cite{CH}.
}
\end{rem}
%%%%%%%%%%%%%%%%%%%%%%%%%%%%%%
%%%%%%%%%%%%%%%%%%%%%%%%%%%%%%
%
%
%
%
%%%%%%%%%%%%%%%%%%%%%%%%%%%%%%%%%
\begin{figure}
\begin{center}
\includegraphics[scale=0.4]{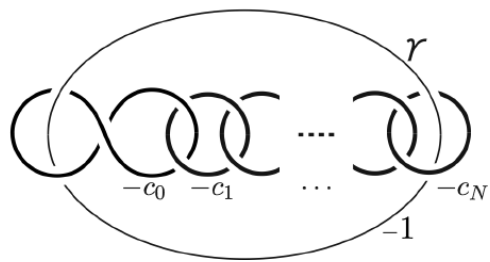} \quad
\caption{\, $C_{p,q}$ with $\gamma$}
\label{fig:total}
\end{center}
\end{figure}
%%%%%%%%%%%%%%%%%%%%%%%%%%%%%%%%%
%
%

%%%%%%%%%%%%%%%%%%%%%%%%%%%%%%%%%%%%%%%%%%%%%%%%%%%%%%%%
%%%                   References
%%%%%%%%%%%%%%%%%%%%%%%%%%%%%%%%%%%%%%%%%%%%%%%%%%%%%%%%
%%%%%%%%%%%%%%%%%%%%%%%%%%%%%%%%%%%%%%%%%%%%%%%%%%%%%%%%

{\small
\par
YAMADA Yuichi \par
Dept. of Systems Engineering, 
The Univ. of Electro-Communications \par
1-5-1,Chofugaoka, Chofu, Tokyo, 182-8585, JAPAN \par
}
{\tt yyyamada@sugaku.e-one.uec.ac.jp} \par

\end{document}